\documentclass[letterpaper, 10 pt, conference]{ieeeconf}  

\IEEEoverridecommandlockouts                              
\usepackage{slashbox}
\usepackage{stfloats}
\usepackage{mathbbold}
\usepackage{amssymb}
\usepackage{graphicx}
\usepackage[version=4]{mhchem}
\usepackage{hyperref}
\usepackage[nameinlink]{cleveref}
\usepackage{xcolor}
\usepackage[ruled,vlined,linesnumbered]{algorithm2e}
\usepackage{upgreek}
\usepackage{bigints}

\def\dd{\text{d}}

\Crefname{algocf}{Algorithm}{Algorithms}
\crefname{algocf}{algorithm}{algorithms}

\title{\LARGE \bf
Effective filtering approach for joint parameter-state estimation in SDEs via Rao-Blackwellization and modularization
}

\author{Zhou Fang, Ankit Gupta, and Mustafa Khammash
\thanks{This work was funded  by the Swiss National Science Foundation under grant number 182653.}
\thanks{Zhou Fang, Ankit Gupta, and Mustafa Khammash are with the Department of Biosystems Science and Engineering at ETH-Zürich.}%
\thanks{\texttt{zhou.fang@bsse.ethz.ch}, \texttt{ankit.gupta@bsse.ethz.ch}, \texttt{mustafa.khammash@bsse.ethz.ch}}
}

\begin{document}

\maketitle
\thispagestyle{empty}
\pagestyle{empty}

\begin{abstract}
Stochastic filtering is a vibrant area of research in both control theory and statistics, with broad applications in many scientific fields. 
Despite its extensive historical development, there still lacks an effective method for joint parameter-state estimation in SDEs.
The state-of-the-art particle filtering methods suffer from either sample degeneracy or information loss, with both issues stemming from the dynamics of the particles generated to represent system parameters.

This paper provides a novel and effective approach for joint parameter-state estimation in SDEs via Rao-Blackwellization and modularization.  
Our method operates in two layers: the first layer estimates the system states using a bootstrap particle filter, and the second layer marginalizes out system parameters explicitly. 
This strategy circumvents the need to generate particles representing system parameters, thereby mitigating their associated problems of sample degeneracy and information loss. 
Moreover, our method employs a modularization approach when integrating out the parameters, which significantly reduces the computational complexity.
All these designs ensure the superior performance of our method.
Finally, a numerical example is presented to illustrate that our method outperforms existing approaches by a large margin.
\end{abstract}

\section{INTRODUCTION}

The advances in measurement technology (e.g., radars \cite{moreira2013tutorial,sun2020mimo}, high-speed cameras \cite{xing2017high}, and microscopes \cite{stephens2003light}) have afforded researchers unprecedented ability to monitor dynamical systems across various scales, from the global climate system to individual living cells. 
However, due to technological and scientific limitations, these technologies cannot directly measure all the dynamical states and system parameters. 
This challenge brings in stochastic filtering, which aims to infer these hidden variables by computing their conditional distributions from noisy partial observations.

Over the past few decades, many effective methods have been developed to address the filtering problem. 
In linear and Gaussian scenarios, the filtering problem can be explicitly solved by the Kalman filter \cite{kalman1960new}. 
For continuous-time Markov chains (CTMCs) with finite state spaces, the conditional distribution can also be explicitly computed by a system of finitely many equations \cite{wonham1964some,confortola2013filtering,chen2022explicitly}. 
In the general setting, the filtering problem is often infinite-dimensional and not explicitly solvable. 
To tackle this challenge, some variants of the Kalman filter (e.g., the extended Kalman filter \cite{jazwinski2007stochastic} and ensemble Kalman filter \cite{evensen1994sequential}) were introduced for nonlinear systems described by stochastic differential equations (SDEs).
For CTMCs having infinite state spaces, the literature \cite{d2022filtered} proposed a method called the filtered finite state projection (FFSP), which approximates the solution by solving the  filtering problem on a truncated state space. 

The FFSP and those variants of the Kalman filter all present certain drawbacks: the FFSP is computationally demanding for large systems, and those Kalman filter variants carry inevitable estimation biases when dealing with non-linear systems.
An alternative to these methods is the bootstrap particle filter (BPF) \cite{gordon1993novel}.
This method employs the simulation-based Monte-Carlo technique to recursively approximates conditional distributions, and it is guaranteed to converge to the exact solution of the filtering problem \cite{crisan2001particle}.
This approach has been successfully used in many applications, including, but not limited to,
wireless communications \cite{djuric2002applications} and biological studies \cite{fang2020stochastic,fang2022stochastic,rathinam2021state}.

Despite its successes, the BPF performs poorly in estimating static variables, e.g., system parameters.
The problem is attributed to the resampling step, which repeatedly reduces the number of distinct particles representing static variables.
Since the BPF does not increase this number in other steps, these particles (representing static variables) soon become identical after a few initial iterations. 
This sample degeneracy diminishes the effectiveness of the Monte-Carlo technique (which relies on a population of distinct particles) and can provide highly inaccurate estimates \cite{doucet2009tutorial}. 

To mitigate sample degeneracy, researchers have developed some improved methods which introduce artificial noise to perturb the particles, thereby increasing the number of distinct particles. 
Such methods include the resample-move method \cite{berzuini2001resample,gilks2001following} and regularized particle filter (RPF) \cite{liu2001combined,oudjane2000progressive,fang2023convergence}.
For joint parameter-state estimation, Crisan and M\'iguez proposed the nested particle filter (NPF) \cite{crisan2018nested}.
This method operates in two layers: the first layer estimates the system parameters using an RPF, and the second layer infers the state variables using a BPF.
Some improvements to this nested approach were reported in  \cite{perez2018probabilistic,perez2021nested}, where the algorithms apply more effective filters in both layers.
Nevertheless, the use of artificial noise (in all the above methods) "throws away" information about parameters, as it assumes parameters (static variables) to be time-varying \cite{liu2001combined}. 
When the particle size is finite, this information loss can also result in highly inaccurate estimates, especially if the artificial noise intensity is not properly chosen.
Consequently, introducing artificial noise might not be the optimal solution to tackle sample degeneracy.

Both sample degeneracy and information loss stem from the particles representing static variables.
From this perspective, a method that avoids generating such particles can effectively tackle these issues.
This strategy aligns with the Rao-Blackwellized particle filter (RB-PF) \cite{murphy2001rao,andrieu2001rao}, which reduces the filter's dependence on Monte Carlo techniques by integrating out some state variables using a finite-dimensional filter.
When the integrated-out variables include the static variables, this method has the potential to effectively mitigate sample degeneracy and information loss. 
The RB-PF also has a nested structure: the first layer employs a BPF to estimate a subset of the hidden variables, and the second layer uses a finite-dimensional filter (e.g. Kalman filter) rather than a particle filter to integrate out the remaining variables. 
The RB-PF has been theoretically shown to have superior performance compared to the BPF in terms of asymptotic variance under quite general conditions \cite{chopin2004central}. 
The RB-PF has been successfully applied to joint parameter-state estimation in biochemical reaction systems described by CTMCs \cite{zechner2014scalable,fang2022scalable}. 
The method in \cite{zechner2014scalable} explicitly marginalizes out all the system parameters using gamma distributions. 
In contrast, the method in \cite{fang2022scalable} employs the FFSP to integrate out all the parameters as well as some state variables, thereby achieving better performance.

Despite its success with CTMCs, the idea of Rao-Blackwellization has not been sufficiently explored for the joint state-parameter estimation in SDEs.
This paper is devoted to filling this gap by providing a principled method  based on Rao-Blackwellization and modularization.
Essentially, our approach adopts the strategy of integrating out system parameters using a finite-dimensional filter.
Specifically, our method uses a nested structure similar to NPF but with the layers reversed.
The first layer of our method infers the system states by employing a BPF.
The second layer applies the  method to the Zakai equation that characterizes the conditional distributions of the system parameters given the state trajectory.
This strategy circumvents the need to generate particles representing system parameters, thereby effectively mitigating the problems associated with sample degeneracy and artificial noise.
Moreover, our method employs a modularization approach for computing the Zakai equation, which significantly reduces the computational effort.
These designs result in an effective method for joint parameter-state estimation in SDEs.

The rest of the paper is organized as follows. 
\Cref{section preliminary} first introduces the mathematical problem of joint parameter-state estimation in SDEs; then it briefly reviews the classical filtering methods for this problem.
In \Cref{section RB-PF}, we introduce our novel approach to this filtering problem.
A numerical example is presented in \Cref{section numerical example} to illustrate the efficiency and accuracy of our method. 
Finally, \Cref{section conclusion} concludes this paper. 
Some terminologies mentioned in this paper are concluded in \Cref{tab:notation}.

\begin{table}[h]
	\centering
		\caption{}
	\begin{tabular}{l l }
		\hline
		Terminology & Meaning  \\
		\hline
		SDE: & Stochastic differential equation \\
		CTMC: & Continuous-time Markov chain \\
		FFSP: & Filtered finite state projection \cite{d2022filtered} \\
		BPF: & Bootstrap particle filter\\
		RPF: & Regularized particle filter\\
		NPF: & Nested particle filter\\
		RB-PF: & Rao-Blackwellized particle filter \\
		$\mathbbold{1}(\cdot)$ : & Indicator function\\
		\hline
	\end{tabular}
	\label{tab:notation}
\end{table}

\section{Stochastic filtering for stochastic differential equations}\label{section preliminary}
We consider stochastic differential systems expressed as:
\begin{align}\label{eq sde}
	\dd X_i(t) =  f_i(\Theta_i, X(t)) \dd t + \sigma_i \dd B_i(t), && i=1,\dots, n
\end{align}
where $X(t)$ is the $n$-dimensional state vector,  $X_i(t)$ is its $i$-th component, $\{f_1, \dots, f_n\}$ are measurable functions, $\{\Theta_1, \dots, \Theta_n\}$ are $n$ unknown system parameters, $\{B_1(t), \dots, B_n(t)\}$ are independent standard Brownian motions, and $\{\sigma_1,\dots, \sigma_n\}$ represent the noise intensities which we assume to be known. 
In many real-word problems, such a system will be measured at consecutive time points $\{t_1, t_2, \dots\}$ with corresponding measurements $Y(t_k)$ ($k=1, 2, \dots$).
We assume that these measurements satisfy
\begin{align*}
	Y(t_k) = h(X(t_k)) + \Sigma W(t_k), && k=1, 2, \dots, 
\end{align*}
where $h(\cdot)$ is a vector-valued measurable function, $\{W(t_1),W(t_2), \dots\}$ are vectors of independent standard Gaussian random variables, and $\Sigma$ is the covariance matrix of the observation noise. 

In practical systems, not all the state variables are measured due to the sensor limitations, which poses a big challenge for better investigation and control of the dynamical system. 
To address this problem, researchers need to infer these hidden states and parameters from the partial observations in real time. 
This mathematical problem, known as stochastic filtering, specifically aims to compute the conditional distribution $\pi_{t_k}(\dd \theta, \dd x) \triangleq\mathbb P \left( \Theta \in \dd \theta, X(t_k) \in \dd x \big| Y(t_1), \dots Y(t_k)\right)$ for $k=1,2, \dots$.
Let us denote the initial distribution as $\pi_{t_0}(\dd \theta, \dd x) \triangleq \mathbb P \left( \Theta \in \dd \theta, X(0) \in \dd x \right)$ and define another conditional distribution $\rho_{t_{k+1}}(\dd \theta, \dd x)\triangleq \mathbb P \left( \Theta \in \dd \theta, X(t_{k+1}) \in \dd x \big| Y(t_1), \dots Y(t_k)\right)$.
Then by Bayes' rule, the target distribution $\pi_{t_k}(\dd \theta, \dd x)$ satisfies the following recursive formulas \cite{crisan2001particle}: 
\begin{align}
	&\rho_{t_{k+1}}(\dd \theta, \dd x) \label{eq. prediction}\\
	&= \sum_{x'}\mathbb P\left( X(t_{k+1}) \in \dd x | \Theta =\theta, X(t_k) = x'\right) \pi_{t_k}(\dd \theta, \dd x') \notag\\
	&\pi_{t_{k+1}}(\dd \theta, \dd x) \propto L\left(Y(t_{k+1}) \big|x\right) \rho_{t_{k+1}}(\dd \theta, \dd x) \label{eq. correction} 
\end{align}
for $k=0,1, \dots$, where the initial time $t_0$ equals to zero, 
and $L(y|x)$ is the likelihood function for the observation given the system state.
The initial distribution $\pi_{t_0}(\dd \theta, \dd x)$ is usually set to be a uniform distribution to reflect our limited knowledge about the specific values of $\Theta$ and $X(0)$. 
In these formulas, we can interpret the formula \eqref{eq. prediction} as the prediction step, which forecasts $X(t_{k+1})$ and $\Theta$ using the observations up to time $t_{k}$; the formula \eqref{eq. correction} can be interpreted as the correction step, which adjusts the predicted distribution based on the new measurement collected at time $t_{k+1}$.

\subsection{Existing filtering methods}

The recursive formulas \eqref{eq. prediction} and \eqref{eq. correction} often cannot be solved explicitly in practical systems because \eqref{eq. prediction} requires the value of the transition probability, which is usually intractable.
This fact necessitates the development of numerical methods for this filtering problem. 
So far, many particle filtering methods have been proposed based on the idea of Monte Carlo. 
We list some as follows. 

\subsubsection{Bootstrap particle filter (BPF) \cite{gordon1993novel}}
The bootstrap particle filter (BPF), also known as sequential importance resampling particle filter, solves \eqref{eq. prediction} and \eqref{eq. correction} by Monte-Carlo samples together with a resampling scheme (see \Cref{alg particle filter} for the detailed algorithm). 
Initially, the algorithm samples $N$ particles from the initial distribution (see Line 1, \Cref{alg particle filter}).
For each $t_k$, the algorithm simulates the particles from time $t_k$ to $t_{k+1}$ according to the dynamical equation (Line 3, \Cref{alg particle filter}).
Then, the empirical distribution of the particles $\texttt{x}_1(t_{k+1}), \dots, \texttt{x}_N(t_{k+1})$ becomes an approximation of the prediction distribution $\rho_{t_{k+1}}(\cdot,\cdot)$ in \eqref{eq. prediction}.
Next, the BPF computes particle weights according to the measurement $Y(t_{k+1})$ (Line 4, \Cref{alg particle filter}) and uses the empirical distribution of the weighted particles to approximate the conditional distribution $\pi_{t_{k+1}}(\cdot,\cdot)$ (Line 5, \Cref{alg particle filter}).
Finally, the algorithm resamples particles to remove non-important samples and reproduce important particles so that the computational complexity is reduced \cite{doucet2009tutorial}. 
The BPF has good reliability in the limit of large particles.
It has been shown that the BPF converges to the exact solution of the filtering problem as the particle size $N$ goes to infinite \cite{crisan2001particle}.

\begin{algorithm}
	\SetAlgoLined
	Sample $N$ particles $\left( \uptheta_1, \texttt{x}_1(0)\right),\dots, \left( \uptheta_N, \texttt{x}_N(0)\right)$ from the initial distribution $\pi_{t_0}$\;
	\For{each time point $t_k$ ($k\in\mathbb Z_{\geq0}$)}{
		Simulate each $\left( \uptheta_j, \texttt{x}_j(\cdot)\right)$ from  $t_k$ to $t_{k+1}$ by \eqref{eq sde}\;
		Compute weights $\texttt{w}_{j}=  L\left(Y(t_{k+1})\big| \texttt{x}_j(t_{k+1})\right)$\;
		Approximated filter: \qquad\qquad\qquad~~ $\bar \pi_{t_{k+1}}(\theta, x)=  \frac{\sum_{j=1}^{N} \texttt{w}_j \mathbbold{1} (\uptheta_j = \theta, \texttt{x}_j(t_{k+1})=x)}{\sum_{j} w_j}$\;
		Resample $\left\{\texttt{w}_{j}, \left( \uptheta_j, \texttt{x}_j(\cdot)\right)\right\}_{j=1,\dots,N}$ to obtain $N$ equally weighted particles, and replace the old particles with these new ones\;
	}
	\caption{Bootstrap particle filter (\cite{gordon1993novel, doucet2009tutorial})}
	\label{alg particle filter}
\end{algorithm}

\subsubsection{Regularized particle filter (RPF)}
Though the BPF is convergent as $N\to \infty$, its performance in estimating system parameters is poor with a finite number of particles due to sample degeneracy.
Specifically, the resampling step reduces the number of distinct particles $\uptheta_1, \dots, \uptheta_N$ used for estimating parameters. 
After several iterations, the BPF often ends up with particles sharing the same parameter part $\uptheta_j$.
This shared  $\uptheta_j$ does not necessarily equal the true parameter values.
Often, random effects can cause this shared $\uptheta_j$ to deviate significantly from the true parameter values, thereby greatly affecting the accuracy of the method.
A more effective alternative to the BPF is the regularized particle filter (RPF) \cite{liu2001combined,west1993approximating,fang2023convergence, oudjane2000progressive}, which introduces some artificial noise to $\uptheta_j$ (in each iteration in the BPF) to ensure greater diversity among the particles $\uptheta_1, \dots, \uptheta_N$. 
The RPF has demonstrated excellent performance in numerous applications (as shown in the aforementioned references), and it also converges to the exact filtering result as $N\to\infty$ under some mild conditions, \cite{fang2023convergence,crisan2014particle,del2000branching, le2004stability}. 

\subsubsection{Nested particle filter (NPF)}
Crisan and M\'{i}guez \cite{crisan2018nested} proposed a two-layer particle filtering algorithm, called the nested particle filter (NPF), for joint estimation of state variables and parameters. The first layer of the NPF employs a regularized particle filter (RPF) to infer the system parameters, i.e., targeting $\pi_{t_k}( \dd \theta)\triangleq \mathbb P \left( \Theta \in \dd \theta \big| Y(t_1), \dots Y(t_k)\right)$.
The second layer uses the BPF to estimate the state variables given a fixed $\theta$, i.e., aiming at $\pi_{t_{k}}(\dd x | \theta) \triangleq \mathbb P \left( X(t_{k})\in \dd x \big| Y(t_1), \dots Y(t_k), \Theta = \theta \right)$.
The two layers are wisely integrated, enabling the NPF to operate in a recursive manner.
Finally, the NPF gives an approximated solution to the filtering problem by combining the  results in both layers according to the Bayes' rule $\pi_{t_k}(\dd \theta, \dd x) = \pi_{t_k}(\dd x|\theta) \pi_{t_k}(\dd \theta)$.
The detailed algorithm for the NPF is provided in \Cref{alg NPF}.
Its validity in the limit of large particles has been shown in \cite{crisan2018nested}.

\begin{algorithm}
	\SetAlgoLined
	Sample $N$ particles $ \uptheta_1,\dots, \uptheta_N$ from $\pi_{t_0} (\dd \theta)$\;
	For each $\uptheta_j$, sample $M$ particles $\texttt{x}^1_j(0), \dots, \texttt{x}^M_j(0)$ from the conditional distribution $\pi_{0}(\dd x | \uptheta_j)$\;
	\For{each time point $t_k$ ($k\in\mathbb Z_{\geq0}$)}{
		\For{each $\uptheta_j$ ($j=1,\dots, N$)}{
			\tcp{Second layer filter}
			Simulate each $ \texttt{x}^\ell_j(\cdot)$ from  $t_k$ to $t_{k+1}$ by \eqref{eq sde} with parameter $\uptheta_j$\;
			Compute weights $\texttt{w}^{\ell}_{j}=  L\left(Y(t_{k+1})\big| \texttt{x}^{\ell}_j(t_{k+1})\right)$\;
			Second-layer filter: $\bar \pi_{t_{k+1}}(x | \uptheta_j) = \frac{\sum_{\ell=1}^M w^{\ell}_j\mathbbold 1(\texttt{x}^{\ell}_j(t_{k+1})=x)}{\sum_{\ell =1}^M w^{\ell}_j}$\;
			Resample $\{w^{\ell}_j, \texttt{x}^{\ell}_j\}_{\ell=1,\dots, M}$ to obtain M equally weighted particles, and replace the old particles with these new ones.
		}
		\tcp{First layer filter}
		Compute weight for each $\theta_j$: $w_j=\sum_{\ell =1}^M w^{\ell}_j$\;
		First-layer filter: $\bar \pi_{t_{k+1}}(\theta) = \frac{\sum_{j=1}^N w_j \mathbbold 1(\uptheta_j =\theta)}{\sum_{j=1}^N w_j}$\;
		Whole filter: $\bar \pi_{t_{k+1}}(\theta, x)=  \bar \pi_{t_{k+1}}(x|\theta) \bar \pi_{t_{k+1}}(\theta)$\;
		Resample $\left\{\texttt{w}_{j}, \left( \uptheta_j, \{\texttt{x}^{\ell}_j(\cdot)\}_\ell\right)\right\}_{j=1,\dots,N}$ to obtain $N$ equally weighted particles, and replace the old particles with these new ones\;
		Perturb each $\uptheta_j$ by some artificial noise
	}
	\caption{Nested particle filter (\cite{crisan2018nested})}
	\label{alg NPF}
\end{algorithm}

\subsubsection{Further remarks on the RPF and NPF}{\label{section further remarks on RPF and NPF}}
In many applications, the performance of RPFs and NPFs largely depends on the wise choice of the artificial noise intensity, which is not easy to determine in advance.
Weak artificial noise cannot effectively circumvent sample degeneracy. 
On the other hand, strong artificial noise could result in the filter ``throwing away" too much information contained in the particles \cite{liu2001combined}. 
Our previous study \cite{fang2023convergence} demonstrates that in some examples, the RPF requires a training step in advance to find the optimal noise intensity, which can be extremely time consuming. 
In summary, employing artificial noise to perturb particles may not be the most effective approach for the joint estimation of state variables and system parameters.

\section{Filtering approach based on Rao-Blackwellization and modularization}\label{section RB-PF}

We propose a new filtering method for joint parameter-state estimation based on Rao-Blackwellization and modularization. 
The key idea is to marginalize out the parameters $\Theta$ using an efficient finite-dimensional filter.
This strategy circumvents the need for generating particles representing $\Theta$, thereby mitigating the related issues.
Moreover, when integrating out $\Theta$, our method employs a divide-and-conquer approach to reduce the computational complexity.
More details about our method are illustrated as follows. 

\subsection{Derivation of our Rao-Blackwellized particle filter} 

First, we give a new formula to re-express the filter $\pi_{t_k}$.
This reformulation integrates out the parameters $\Theta$ following the idea of Rao-Blackwellization \cite{rao1992information,blackwell1947conditional}.
Let $X_{0:t}$ be the whole trajectory of $X(\cdot)$ from time zero to $t$, $Y_{t_1:t_k}$ the measurements up to $t_k$, and $\pi_{\Theta|X} (t, \dd \theta) = \mathbb P \left(
\Theta \in \dd \theta \big| X_{0:t}
\right)$.
Then, by the tower property, we can express $\pi_{t_k}$ by
\begin{align}\label{eq. pi RB}
	&\pi_{t_k} (\dd \theta, \dd x)\notag\\
	&= \mathbb E\left[
	          \mathbb P \left(
	              \Theta \in \dd \theta \big| X_{t_0:t_k}, Y_{t_1: t_k} 
	          \right)\,
	          \mathbbold 1\left( X(t_{k}) \in \dd x\right)
	          \Big| Y_{t_1:t_k} 
	\right]\notag\\
	& = \mathbb E\left[
	\pi_{\Theta|X}(t_k,\dd \theta)\,
	\mathbbold 1\left( X(t_{k}) \in \dd x\right)
	\Big| Y_{t_1:t_k} 
	\right] 
\end{align}
where the second equality holds because the likelihood function $L(y|x)$ does not depend on the parameters in our problem. 
Similarly, the $\rho_{t_{k+1}}$ can be expressed by 
\begin{align}\label{eq rho RB}
	&\rho_{t_{k+1}} (\dd \theta, \dd x) \notag\\
	& = \mathbb E\left[
	\pi_{\Theta|X}(t_{k+1},\dd \theta)\,
	\mathbbold 1\left( X(t_{k+1}) \in \dd x\right)
	\Big| Y_{t_1:t_k} 
	\right].
\end{align}

Equations \eqref{eq. pi RB} and \eqref{eq rho RB} suggest that the filtering problem can be numerically solved by generating samples for the process $\left(\pi_{\Theta|X}(t,\cdot), X(t)\right)$.
This scheme will compute the conditional distribution $\pi_{\Theta|X}(t,\cdot)$  for each simulated trajectory rather than generating the particles that represent parameters $\Theta$.
Consequently, this can mitigate the issues associated with sample degeneracy and artificial noise occurring in the aforementioned particle filters.

Next, we introduce a divide-and-conquer approach to solve $\pi_{\Theta|X}(t,\cdot)$.
When the terms in the dynamical equation \eqref{eq sde} are regular enough, the density of $\pi_{\Theta|X}(t,\cdot)$ is the unique normalized solution of the Zakai equation \cite{bain2009fundamentals}
\begin{align*}
	\left\{
	\begin{array}{l}
			\dd \rho_{\Theta|X}(t, \theta) =  \rho_{\Theta|X}(t, \theta) \sum_{i=1}^n \frac{f_i(\theta_i, X(t))}{\sigma^2_1}\dd X_i(t)\\
			\rho_{\Theta|X}(0, \theta) \dd\theta  =  \pi_{\Theta|X}(0, \dd \theta) 
	\end{array}
		\right.
\end{align*}
where $\theta_i$ is the $i$-th component of $\theta$.
Basically, this Zakai equation is high dimensional, and it suffers the curse of dimensionality when solved directly using grid-based methods. 
Fortunately, our system \eqref{eq sde} has a nice structure where the $i$-th parameter only immediately affects the $i$-th state.
This enables a modularization method for computing the Zakai equation.
Specifically, the Zakai equation suggests 
\begin{align*}
	&\rho_{\Theta|X}(t, \theta)
	\propto
	~\rho_{\Theta|X}(0, \theta)  \times \\
	&\quad~\prod_{i=1}^{n}
	\exp\left\{
			 \int_0^t \frac{f_i(\theta_i, X(t))}{\sigma^2_1} \left[\dd X_i(t)  
			 - \frac{f_i(\theta_i, X(t))}{2} \dd t
			 \right] 
	\right\}.
\end{align*}
Recall that we assume a uniform prior distribution for $\Theta$ and $X(0)$.
This means all the parameters are conditionally independent given $X(0)$, and 
$\rho_{\Theta|X}(0, \theta)$ can be written by $\rho_{\Theta|X}(0, \theta) = \prod_{i=1}^n \rho_{\Theta_i|X}(0, \theta_i) $ where $\rho_{\Theta_i|X}(0, \cdot)$ is the marginal conditional distribution for $\Theta_i$ given $X(0)$.
Thus, the solution of the Zakai equation can be expressed by 
\begin{align}\label{eq. divide and conquer for parameters}
	&\rho_{\Theta|X}(t, \theta) \propto \\
	&
	{\footnotesize
	\prod_{i=1}^{n} \underbrace{\rho_{\Theta_i|X}(0, \theta_i)  
	\exp\left\{
		\int_0^t \frac{f_i(\theta_i, X(t))}{\sigma^2_1} \left[\dd X_i(t)  
						- \frac{f_i(\theta_i, X(t))}{2} \dd t
						\right] 
	\right\}
	}_{=: ~\rho_{\Theta_i|X}(t, \theta)} 
	\notag}
\end{align}
suggesting that the parameters are conditionally independent given $X_{0:t}$. 
Here, we denote $\rho_{\Theta_i|X}(t, \cdot)$ as the un-normalized marginal conditional distribution for $\Theta_i$.
With this conditional independence, we can compute $\pi_{\Theta|X}(t,\cdot)$  (or equivalently, $\rho_{\Theta|X}(t, \cdot)$) by applying the Euler–Maruyama method to each marginal distribution $\rho_{\Theta_i|X}(t, \cdot)$ rather than the joint distribution.
This strategy reduces the computational complexity from $O\left(\texttt{G}^n\right)$ to $O\left(n \texttt{G}\right)$, with $G$ the number of grid points for each parameter. 
Consequently, this divide-and-conquer method, as suggested by \eqref{eq. divide and conquer for parameters}, is efficient even in  high-dimensional cases.

\subsection{Algorithm of our Rao-Blackwellized particle filter}

Following the idea presented above, we provide a Rao-Blackwellized particle filter (RB-PF) in \Cref{alg RB-PF}.
Essentially, we generate samples $\left\{\left(\texttt{x}_j(\cdot), \bar \rho^1_j, \dots, \rho^n_j \right)\right\} _{j=1,\dots, N}$ for the processes $\left(X(\cdot), \rho_{\Theta_1|X}, \dots, \rho_{\Theta_n|X} \right)$ and use them to approximate the exact filter $\pi_{t_k}(\cdot)$ by \eqref{eq. pi RB}.
Further details are elaborated as follows.

\begin{algorithm}[htb]
	\SetAlgoLined
	Sample N particles $\texttt{x}_1(0), \dots, \texttt{x}_N(0)$ from $\pi_{t_0} $ \;
	For every $\Theta_i$, select a finite set $\mathbbold \Theta_i\subset\mathbbold R$. Then, for each $\texttt{x}_j(0)$ and each $i$, denote marginal distribution $\bar \rho^{i}_{j} (0, \theta_i) \propto \frac{  \mathbb P\left( \Theta_i= \dd \theta_i | X(0) = \texttt{x}_j(0)\right)}{\dd \theta_i}$, $\forall \theta_i \in \mathbbold \Theta_i$ \;
	\For{each time point $t_k$ ($k\in\mathbb Z_{\geq0}$)}{
		Simulate each $\texttt{x}_j(\cdot)$ from $t_k$ to $t_{k+1}$ by \eqref{eq sde} with parameters $\uptheta_j$ sampled from the distribution $\bar \rho_{j} (t_k, \theta) \triangleq \prod_{i=1}^n \bar \rho^{i}_{j} (t_k, \theta_i)$\;
		For each $\texttt{x}_j(\cdot)$ and $i\in\{1,\dots,n\}$, compute 
		{\scriptsize
			$$\bar \rho^{i}_{j} (t_{k+1}, \theta_i) = \qquad\qquad\qquad\qquad
			\qquad\qquad\qquad\qquad\qquad\qquad\qquad$$
			$$ \qquad\quad\bar \rho^{i}_{j} (t_k, \theta_i) \text{e}^{
				\bigintssss_{t_k}^{t_{k+1}} \frac{f_i(\theta_i, \texttt{x}_j(s))}{\sigma^2_i} \left[ \dd \left(\texttt{x}_j\right)_i(s)
				- \frac{f_i(\theta_i, \texttt{x}_j(s))}{2} \dd t \right]
			}$$} for every $\theta_i \in \mathbbold \Theta$, and then normalize it\;
		Compute weights $\texttt{w}_{j}=  L\left(Y(t_{k+1})\big| \texttt{x}_j(t_{k+1})\right)$\;
		Approximated filter:
		{\footnotesize
			$$\bar \pi_{t_{k+1}}(\theta, x)=  \frac{\sum_{j=1}^{N} \left[\texttt{w}_j \mathbbold{1} ( \texttt{x}_j(t_{k+1})=x) \prod_{i=1}^n \bar\rho^i_{j}(t_{k+1}, \theta_i)\right]}{\sum_{j} w_j}$$}\\
		Reample $\left\{\texttt{w}_{j}, \left( \texttt{x}_j, \bar \rho^{1}_{j}, \dots, \bar \rho^{n}_{j}  \right)\right\}_{j=1,\dots,N}$ to obtain $N$ equally weighted particles, and replace the old particles with these new ones\;
	}
	\caption{Rao-Blackwellized particle filter}
	\label{alg RB-PF}
\end{algorithm}

First, the algorithm generates particle $\texttt{x}_1(0), \dots, \texttt{x}_N(0)$ from the initial distribution (Line 1, \Cref{alg RB-PF}).
Then, for each $\texttt{x}_j(0)$, the algorithm creates $\bar \rho^i_j(0, \cdot)$ to represent the conditional density of $\Theta_i$ given $X(0)=\texttt{x}_j(0)$ (Line 2, \Cref{alg RB-PF}).
Due to memory constraints, a digital computer cannot store all the values of this conditional density function. 
Consequently, for each $\Theta_i$, the algorithm only stores the values of the density function corresponding to some representative points selected within $\mathbb R$.
The set containing these selected points is denoted by $\mathbbold \Theta_i$.

Then, at each time point $t_k$, \Cref{alg RB-PF} solves the prediction and correct equations \eqref{eq. prediction} and \eqref{eq. correction} by simulating the particles according to \eqref{eq sde} and \eqref{eq. divide and conquer for parameters}.
The algorithm simulates every $\texttt{x}_j(\cdot)$ from time $t_k$ to $t_{k+1}$ according to the dynamics \eqref{eq sde} with parameters sampled from the conditional distribution $\bar \rho_{j} (t_k, \theta) \triangleq \prod_{i=1}^n \bar \rho^{i}_{j} (t_k, \theta_i)$ (Line 4, \Cref{alg RB-PF}).
It can be easily shown that given $X(t_k) = \texttt{x}_j(t_k)$ and  $\bar \rho_{j} (t_k, \cdot)$ exactly equaling  $\rho_{\Theta|X}(t_k, \cdot)$ (up to normalization), the trajectory of $\texttt{x}_j(\cdot)$ from $t_k$ to $t_{k+1}$ (as produced by our algorithm) has the same distribution as $X_{t_k:t_{k+1}}$.
When $\bar \rho_{j} (t_k, \cdot)$ provides an accurate, but not perfect, approximation of $\rho_{\Theta|X}(t_k, \cdot)$, the generated $\texttt{x}_j(\cdot)$ should still statistically resemble $X(\cdot)$. 
After the simulation of $\texttt{x}_j(\cdot)$, the algorithm computes the marginal conditional densities $\bar \rho^{i}_{j} (t_{k+1}, \cdot)$ according to \eqref{eq. divide and conquer for parameters} (Line 5, \Cref{alg RB-PF}).
By \eqref{eq rho RB}, the particles $\left\{\left(\texttt{x}_j(t_{k+1}), \bar \rho^1_j(t_{k+1}, \cdot), \dots, \rho^n_j(t_{k+1}, \cdot) \right)\right\} _{j=1,\dots, N}$ can approximate the predicting distribution $\rho_{t_{k+1}}$ using $\bar \rho_{t_{k+1}}(\theta, x)=  \frac{\sum_{j=1}^{N} \left[ \mathbbold{1} ( \texttt{x}_j(t_{k+1})=x) \prod_{i=1}^n \bar\rho^i_{j}(t_{k+1}, \theta_i)\right]}{N}$.
To approximate the filter $ \pi_{t_{k+1}}(\theta, x)$, our algorithm computes weights for all the particles (Line 6, \Cref{alg RB-PF}) and provides a filter 
$\bar \pi_{t_{k+1}}(\theta, x)$ according to \eqref{eq. correction} and \eqref{eq. pi RB} (Line 7, \Cref{alg RB-PF}).
Finally, our algorithm resamples the particles to accelerate the speed.

\subsection{Some discussions about our method}

Our RB-PF also has a nested structure.
Concretely, the particles $\texttt{x}_1(\cdot), \dots, \texttt{x}_N(\cdots)$ can be seen as the first layer estimating the system state, and $\left\{\left(\bar \rho^1_j, \dots, \rho^n_j \right)\right\} _{j=1,\dots, N}$ form the second layer estimating the parameters (given the state trajectories). 
Compared with the NPF (\Cref{alg NPF}), our filter reversed its order of layers for estimating the parameters and states.
More importantly, the second layer in our algorithm uses a finite-dimensional filter (rather than a BPF or RPF) for parameter estimation.
This strategy circumvents the need to generate particles representing $\Theta$, thereby mitigating issues related to sample degeneracy and artificial noise.

Our algorithm improves classical particle filtering methods (e.g., the BPF, RPF, and NPF) at the cost of requiring more computational resources for the same particle size.
Essentially, our RB-PF needs to additionally compute and store the conditional distributions $\bar \rho^i_j(t, \cdot)$, which necessitates more computational time and computer memory. 
We provided a modularization-based approach for these conditional distributions, thereby reducing the additional computational costs to some extent. 
Still, the extra computational effort is not negligible.
Nonetheless, this additional computation is worthwhile.
It can effectively mitigate the problems associated with RPF and NPF, and, consequently, our method can outperform existing methods for the same computational time. 
The next section will further illustrate this point using a numerical example.

\section{Numerical Example}\label{section example}\label{section numerical example}

Here, we illustrate the superior performance of our method using the stochastic Lorenz-63 model. 
The model consists of three states ($X_1(t),X_2(t), X_3(t)$) and three parameters ($\Theta_1,\Theta_2, \Theta_3$). 
Its dynamics is described by 
\begin{align*}
	\dd X_1(t) &= -\Theta_1 \left[X_1(t) - X_2(t)\right] \dd t + \sigma \dd B_1(t) \\
	\dd X_2(t) &=  \left[\Theta_2 X_1(t) -X_2(t) -X_1(t)X_3(t)\right] \dd t + \sigma \dd B_2(t) \\
	\dd X_3(t) &= \left[X_1(t)X_2(t) - \Theta_3 X_3(t)\right] \dd t + \sigma \dd B_3(t)
\end{align*}
where $\sigma$ is a known parameter depicting the intensity of the process noise, and $B_1(t)$, $B_2(t)$, and $B_3(t)$ are independent standard Brownian motions.
Here, we consider $\sigma =1$.
Clearly, this system conforms to the model \eqref{eq sde}.
We assume that the process are measured at time points $\{0.05, 0.10,\dots, 10\}$, and the measurements $Y(t_k)$ satisfies
\begin{align*}
	Y(t_k) = \left[
	\begin{array}{c}
	X_1(t_k) \\ 
	X_3(t_k)
	\end{array}
	\right] + 
	\left[
	\begin{array}{c}
		W_1(t_k) \\
		W_2(t_k)
	\end{array}
	\right] 
	&& \text{for } t_{k} = 0.05,\dots.
\end{align*}
Here, $	\{W_1(t_k),W_2(t_k)\}_{k=1,\dots, 200}$ are independent standard Gaussian noise. 
The aim of this numerical example is to infer the state and parameters in real time, i.e., to compute $\pi_{t_k}(\dd \theta, \dd x) \triangleq\mathbb P \left( \Theta \in \dd \theta, X(t_k) \in \dd x \big| Y(t_1), \dots Y(t_k)\right)$.

\begin{figure}[htb]
	\centering
	\includegraphics[width=\linewidth]{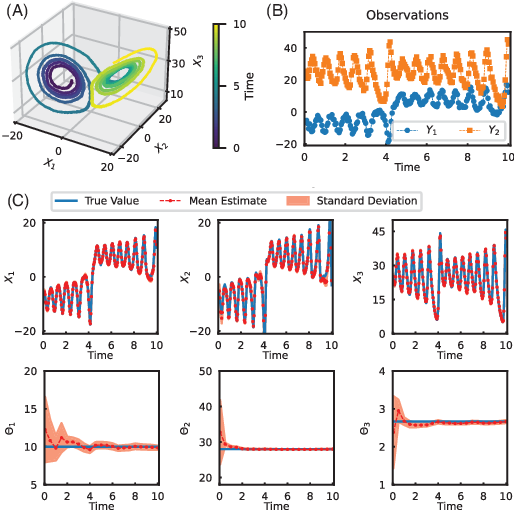}
	\caption{Results of our RB-PF applied to the stochastic Lorenz-63 system.
	(A) Simulated trajectory of the system states. The trajectory begins at the dark-color end and terminates at the light-color end.
	(B) Observation process. 
	(C) Results in estimating state variables and parameters. 
	The blue line denotes the true value, the red line represents the mean estimate provided by our RB-PF, and the light red region indicates the area within one standard deviation. 
	}
	\label{fig RBPF result}
\end{figure}

First, we examined the performance of our RB-PF in solving this filtering problem. 
We assumed that the initial state and parameters were independent and satisfied some uniform distributions: $X_1(0) \sim \mathcal U (-9, -3)$, $X_2(0)\sim \mathcal U (-9, -3)$, $X_3(0)\sim \mathcal U(20, 28)$, $\Theta_1 \sim \mathcal U(5, 20)$, 
$\Theta_2 \sim \mathcal U(18, 50)$, and $\Theta_3 \sim \mathcal{U}(1, 8)$.
Then, we simulated a trajectory of states and observations with the initial state $\left(-6, -5, 24.5\right)^{\top}$ and parameters $\left(10, 28, 8/3\right)$.
Finally, we applied our RB-PF (with a particle size of $N =20,000$) to infer the system state and parameters from these simulated measurements.
The filtering algorithm was performed on the Euler Computing Cluster at ETH Zurich, using a node with 12-core CPUs. 
The whole computational time was approximately 55 minutes.
The numerical results are presented in \Cref{fig RBPF result}.

The numerical result shows that our approach provides sharp estimates for the system states and parameters.
The mean estimates of the state variables almost overlap with the true state trajectories, with the standard deviations too small to be visible (the first row in \Cref{fig RBPF result}.(C)).
The mean estimates for the system parameters also fast converge to the true values (the second row in \Cref{fig RBPF result}.(C)).
From time four onward, these estimates closely match the true parameter value, with very small standard deviations.
All these results demonstrate the accuracy of our approach in solving this filtering problem.

\begin{figure}[htb]
	\centering
	\includegraphics[width=\linewidth]{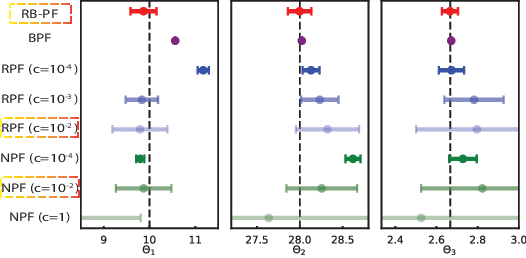}
	\caption{Comparison of different filtering methods in  parameter estimation.
	The black dash line represents the true values of these system parameters. 
	The colored lines present the estimates from different filters, with the dots indicating the mean estimates and the error bars representing the standard deviation.
	The names of the particle filters are listed on the left, where $c$ denotes the intensity of the artificial noise.
	The three highlighted filters are the ones where all the confidence intervals (mean $\pm$ one standard deviation) include their respective true parameter values.
	}
	\label{fig comparison}
\end{figure}

Next, we compared our RB-PF with other competing approaches (the BPF, RPF, and NPF).
To ensure a fair comparison, we applied these filters to the same observation trajectory (as shown in \Cref{fig RBPF result}.B), and we carefully selected their particle sizes so that their computational time was similar to that of the RB-PF (approximately 55 minutes).
Specifically, the sample sizes of the BPF and RPF were set to be 40,000.
For the NPF, the particle sizes in the first layer and second layer ($N$ and $M$, respectively) were set to be equal as suggested in \cite{crisan2018nested}, with the specific value chosen to be 200.
The artificial noise in the RPF and NPF was generated from a normal distribution with mean $(0, 0, 0)^{\top}$ and covariance matrix $c I_{d}$.
Here, $c$ is a tunable hyper-parameter, and $I_d$ is the identity matrix.
To avoid the perturbed particle leaving the defined parameter region $[5,20]\times[18, 50]\times[1, 8]$, we repeatedly generated artificial noise for each particle until the perturbed particle remained within this parameter region.
The comparison results are presented in \Cref{fig comparison}.

The result shows that the RB-PF is the most accurate in estimating system parameters. 
First, the RB-PF is one of the three filters where all the confidence intervals (mean $\pm$ one standard deviation) include their respective true parameter values. 
Moreover, among these three filters, the RB-PF has the smallest standard deviations, suggesting it is the most accurate one. 

\Cref{fig comparison} also reveals the great difficulty of choosing a proper artificial noise intensity ($c$) for the RPF and NPF to balance between sample degeneracy and information loss.
When $c$ is small, these filters still suffer sample degeneracy, resulting in very similar particles whose confidence interval does not necessarily covers the true parameter value (see the RPFs and NPFs with $c < 10^{-2}$).
Particularly, when $c= 0$, the RPF degenerates to the BPF, resulting in particles with identical parameter parts.
As $c$ increases, the artificial noise can lead to significant information loss, causing a large standard deviation in parameter estimation (see the RPFs and NPFs with $c\geq 10^{-2}$).

Thanks to Rao-Blackwellization, our method do not generate particles representing static variables $\Theta$ and, therefore, avoids the problems related to sample degeneracy and artificial noise.
This is one contributing factor to the optimal performance of our approach.
Recall that this improvement comes at the cost of requiring more computational resources for computing the Zakai equation.
By employing the divide-and-conquer strategy (as suggested by \eqref{eq. divide and conquer for parameters}), this additional computation only leads to a manageable increase in computational effort.
In this example, the computational time required by the additional computation is  comparable to that spent on the remaining part of the algorithm.
Consequently, for the same computational time, the RB-PF can still employ half of the sample size compared to the other competing methods, resulting in its superior performance in parameter estimation.
Overall, the additional computational resources required in RB-PF are acceptable and yield significant benefits.

\section{Conclusion}\label{section conclusion}
Facing the challenge of joint parameter-state estimation in SDEs, we proposed a novel and effective filtering approach based on Rao-Blackwellization and modularization. 
Our method operates in two layers: the first layer estimates state variables using a BPF, and the second layer integrates out all the parameters using a Euler–Maruyama method. 
This strategy eliminates the need to generate particles representing parameters and, therefore, circumvents the problems of sample degeneracy and information loss presented in the state-of-the-art methods.
Moreover, our method employs a modularization approach in the second layer, which significantly reduces the required additional computational effort.
These designs result in an effective filtering algorithm for joint parameter-state estimation in SDEs. 
Its superior performance was also demonstrated through a numerical example.

There are a few topics deserving further investigation in future work. 
First, a theoretical analysis of this method is needed to investigate its convergence, asymptotic variance, and limitations.
Second, the method can be further improved by integrating out some state variables in addition to the system parameters. 
Our previous work \cite{fang2022scalable} can be beneficial for this extension, as it successfully employed this idea in the estimation of CTMCs.


\end{document}